\documentclass{article}
\title{{\bf MONOMIAL RESOLUTIONS OF MORPHISMS OF ALGEBRAIC SURFACES}}

\def\hfl#1{\smash{\mathop{\hbox to 12mm{\lefttarrowfill}}
\limits^{\scriptstyle#1}}}
\def\hfr#1{\smash{\mathop{\hbox to 12mm{\rightarrowfill}}
\limits^{\scriptstyle#1}}}
\def\vfd#1#2{\llap{$\scriptstyle #1$}\left\downarrow
\vbox to 6mm{}\right.\rlap{$\scriptstyle #2$}}

\def\diagram#1{\def\normalbaselines{\baselineskip=0pt
\lineskip=10pt\lineskiplimit=1pt} \matrix{#1}}
\begin{document}

\author{Dale Cutkosky\thanks{The first author was partially supported by NSF}\\
Department of Mathematics\\
 University of Missouri\\
Columbia, MO 65211, USA
\and
Olivier Piltant\\
Centre de Math\'ematiques\\
\'Ecole polytechnique\\
F91128 PALAISEAU Cedex\\
France
}

\maketitle
\vskip 5mm

\newtheorem{defin}{Definition}
\newtheorem{prop}{Proposition}
\newtheorem{lem}{Lemma}
\newtheorem{cor}{Corollary}
\newtheorem{them}{Theorem}
\newtheorem{conj}{Conjecture}
\centerline{\it Dedicated to Robin Hartshorne on the occasion of his sixtieth
birthday} \vskip 1cm

\section{Introduction.}
\vskip 0mm
Let $k$ be a perfect field and $L/K$ be a finite separable field extension of
one-dimensional function fields over $k$. A classical result (c.f. I.6, [Ha])
states that
$K$ (resp. $L$) has a unique proper and smooth model $C$ (resp. $D$),
and that there is a unique morphism of curves $f: \ D \rightarrow C$
inducing the field inclusion $K \subset L$ at the generic points of $C$ and
$D$. It has the following properties:

\begin{itemize}
\item[(i)]
$f$ is a finite morphism.
\item[(ii)]
$f$ is {\it monomial} on its tamely ramified locus; let $\beta \in D$ be any
point, with $\alpha:=f(\beta) \in C$, such that the extension of discrete
valuation rings ${\cal O}_{Y,\beta} / {\cal O}_{X,\alpha}$ is tamely
ramified. There exists a local-\'etale ring extension $R$ of
${\cal O}_{Y,\beta}$ and regular parameters $u$ of ${\cal O}_{X,\alpha}$
and $\overline{x}$ of $R$ such that $u=\overline{x}^a$ for some $a$ prime to
the characteristic of $k$.
\end{itemize}

In this paper, we investigate a two-dimensional version of this statement,
that is, $L/K$ is a finite separable field extension of two-dimensional
function fields over $k$. By birational resolution of singularities ([Ab4],
[H], [Li2])
and elimination of indeterminacies (theorem 26.1, [Li]),
there exists a proper and smooth model $X$ (resp. $Y$)
of $K$ (resp. $L$), together with a morphism $f: \ Y \rightarrow X$ inducing
the field inclusion $K \subset L$ at the generic points of $X$ and $Y$.
Such a morphism is in general neither finite nor monomial. \\

In [Ab2], the question is raised of whether this can be arranged by
blowing-up: does there exist compositions of point blow-ups $Y' \rightarrow Y$
and $X' \rightarrow X$, together with a map $f': \ Y' \rightarrow X'$ such
that $f'$ is finite and/or monomial? It is actually shown (theorem 12, [Ab2])
that $f'$ finite cannot in general be achieved. The obstruction is local
for the Riemann-Zariski manifold of $L/k$.\\

This leaves open the question of whether $f'$ can be
taken to be monomial. For complex surfaces, a positive answer has been given
in [AKi] (theorem 7.4.1). Their method however does not generalize to
positive characteristic, due to the lack of canonical forms for the
equations defining $f$ (assertion 7.4.1.1, [AKi]).\\

We present a quite general solution to this problem: any {\it proper, tamely
ramified} morphism $f: \ Y \rightarrow X$ of surfaces (which are separated but
not necessarily proper over $k$), inducing the field inclusion $K \subset L$
at the generic points of $X$ and $Y$, can, after performing suitable
compositions of point blow-ups $Y' \rightarrow Y$ and $X' \rightarrow X$,
be arranged to a monomial morphism $f': \ Y' \rightarrow X'$.
Moreover, there is a unique {\it minimal} such $f'$. \\

Our method is constructive. That is, we give an algorithm, which, starting
from an arbitrary proper $f$ as above, produces its associated minimal $f'$.
This algorithm is explained in section 4.
An easy reduction (proposition 8) shows that it can be assumed
that both of the critical locus $C_f$ and the branch locus
of $f$ are divisors with strict normal crossings.\\
In section 3, we then attach to every vertical component $E$ of $C_f$
a nonnegative integer,
its complexity $i_E$ (definition 4), which is zero if $f$
is monomial at all points of $E$ (compare with that used in [AKi], p.222).
That our algorithm eventually makes it drop, which is the main technical point,
follows from propositions 6 and 7. \\
The main theorem is stated in section 2,
together with the appropriate notions of monomial and of tamely ramified
(not necessarily finite) morphism.\\

In Section 5, we give a proof that $f$ can be made toroidal when $k$ is an
algebraically closed field of
 characteristic zero.
Although this result is known, for instance it is implicit in [AKi], we include
it as
an interesting point in the general theory of resolution of morphisms of
surfaces.\\

There is a local formulation of a monomial resolution for a mapping.
Suppose that $f:Y\rightarrow X$ is a morphism of varieties over a field $k$. If
$f(p)=q$, we have an induced
homomorphism of local rings
$$
R={\cal O}_{X,q}\subset S={\cal O}_{Y,p}
$$
We will say that $R\rightarrow S$ is a monomial mapping if there are regular
parameters $(x_1, .... ,x_m)$
in $R$,  $(y_1, ... ,y_n)$ in $S$ (with $m\le n$), units
$\delta_1,\ldots,\delta_n\in S'$ and a matrix $(a_{ij})$ of
nonnegative integers such that  $(a_{ij})$ has  rank $m$, and
$$
\begin{array}{lll}
x_1 &=& y_1^{a_{11}} ..... y_n^{a_{1n}}\delta_1\\
&&\vdots\\
x_m &=& y_1^{a_{m1}} ..... y_n^{a_{mn}}\delta_m.
\end{array}\leqno (0.1)
$$
Suppose that $V$ is a valuation ring of the quotient field $K$ of $S$, such
that $V$ dominates $S$.
Then we can ask if there are sequences of monoidal transforms $R\rightarrow R'$
and $S\rightarrow S'$ such that
$V$ dominates $S'$, $S'$ dominates $R'$, and $R'\rightarrow S'$ has an
especially good form.
$$
\begin{array}{lll}
R'&\rightarrow&S'\subset V\\
\uparrow&&\uparrow\\
R&\rightarrow &S
\end{array}\leqno (0.2)
$$
Zariski's Local Uniformization Theorem [Z1] says that (when $\mbox{char}(k)=0$)
there exists a diagram (0.2)
such that $R'$ and $S'$ are regular.

In Theorem 1.1 [C]  we  obtain a diagram (0.2) making $R'\rightarrow S'$ a
monomial mapping whenever the quotient field of $S$ is a finite extension of
the quotient field of $R$, and
the characteristic of $k$ is 0.

If  $R'\rightarrow S'$ is a
 mapping of the form (0.1), and the characteristic of $k$ is zero, there exists
a local etale extension $S'\rightarrow S''$
such that $S''$ has regular parameters $\overline y_1,\ldots \overline y_n$
such that
$$
\begin{array}{lll}
x_1 &=& \overline y_1^{a_{11}} ..... \overline y_n^{a_{1n}}\\
&&\vdots\\
x_m &=& \overline y_1^{a_{m1}} ..... \overline y_n^{a_{mn}}.
\end{array}\leqno (0.3)
$$
In char $p>0$, the form (0.3) is not possible to obtain from a monomial mapping
by an \'etale extension in general.
Already in dimension 1,
$$
x=y^p+y^{p+1}
$$
gives a simple counterexample. However, the above example is a monomial
mapping. In fact, if $R$ and $S$
 are regular local rings of dimension 1, then
$R\subset S$ is a monomial mapping, since $R$ and $S$ are Dedekind domains.

If $R$ and $S$ have dimension 2,  $k$ is a field of characteristic $p>0$,and
$V$ is a valuation ring
 dominating $S$, then we ask if
it is  possible to obtain a diagram (0.2)  making $R'\rightarrow S'$ a
monomial mapping.  From our theorem 1, we deduce a positive answer whenever $p$
does not divide the
order of a Galois closure of the quotient field of $S$ over the quotient field
of $R$.

\section{Preliminaries and statement of main result.}
\vskip 0mm
All along this article, $k$ denotes a perfect field of characteristic
$p \geq 0$, and $K / k$ a finitely generated field extension. $L / K$ is
a finite separable field extension.

By an algebraic $k$-scheme, we mean a Noetherian separated $k$-scheme, all
whose local rings are essentially of finite type over $k$. The function
field of an integral algebraic $k$-scheme $X$ is denoted by $K(X)$. If
$\alpha$ is a closed point of such a scheme, its ideal sheaf is denoted by
$M_\alpha$.  If $R$ is a local ring, its residue field is denoted
by $\kappa (R)$.

\begin{defin}
A proper, generically finite morphism of integral algebraic $k$-schemes
$f: \ Y \rightarrow X$ is called a model of the field extension $L/K$ if
$K(X)=K$, $K(Y)=L$, ${\rm dim}X \ (= {\rm dim}Y) = {\rm tr.deg}_kK$, and if
the following diagram commutes:
$$
\diagram{ Spec L& \hfr{}{}&Spec K \cr
\vfd{}{} & & \vfd{}{} \cr
Y & \hfr{f}{}& X \cr}
$$
\end{defin}

A model is said to be proper if $X /k$ (and hence $Y/k$ as well) is proper,
and nonsingular if both of $X$ and $Y$ are nonsingular. Models are partially
ordered by domination, where a model $f': \ Y' \rightarrow X'$ dominates
another model $f: \ Y \rightarrow X$ if there exist {\it proper} maps
$\pi: \ Y' \rightarrow Y$ and $\eta: \ X' \rightarrow X$ such that the
following diagram commutes (and is compatible with the maps of definition 1):
$$
\diagram{ Y'& \hfr{f'}{}& X' \cr
\vfd{\pi}{} & & \vfd{\eta}{} \cr
Y & \hfr{f}{}& X \cr} 
$$

Recall that a $k$-valuation ring $V$ ($k \subset V$) of $L$, with
$V \subset L=K(V)$ is said to be {\it divisorial} if its group is
isomorphic to ${\bf Z}$, and if ${\rm tr.deg}_k \kappa(V)  = {\rm tr.deg}_k
L -1$ (divisorial valuations are called {\it prime divisors} in [ZS2], p.88). 
A generically finite inclusion $W \subset V$ of divisorial
$k$-valuation rings is said to be tamely ramified if ${\rm char}k=0$, or
if ${\rm char}k=p>0$, its ramification index is not divisible by $p$,
and the residue field extension $\kappa(V) / \kappa (W)$ is separable.

\begin{defin}
A model $f: Y \rightarrow X$ is said to be tamely ramified if for
every divisorial $k$-valuation
ring $V$ of $L$, with $K(V)=L$, having a center in $Y$, the extension of
valuation rings $V/ V \cap K$ is tamely ramified.
\end{defin}

{\it Remark:} since the models we are considering are not necessarily finite, 
a notion of tame ramification involving {\it all} divisorial valuations rings 
having a center in $Y$ is needed. For finite morphisms, the usual definition 
(2.2.2 of [GM], or p.41 of [Mi]) only involves those divisorial valuations 
as above whose center in $X$ has codimension one.\\
This raises the following problem: if $f: \ Y \rightarrow X$ is a model, and 
if the induced finite map
\[
\overline{f}: \ {\bf Spec}f_*{\cal O}_Y \rightarrow X
\]
is tamely ramified in the sense of [GM], under which conditions is it true 
that $f$ is tamely ramified according to definition 2?\\

From now on, it will be assumed that $tr.deg_k K =2$. All models therefore
are proper, generically finite morphisms of integral surfaces.\\

Given a {\it nonsingular} model $f: \ Y \rightarrow X$, its critical
locus is denoted by $C_f$.
A scheme structure on $C_f$ is given by the vanishing of the Jacobian
determinant. Since $L/K$ is separable, $C_f$ is a divisor on $Y$.
There exist well defined effective divisors $R_f, S_f$ on $Y$ such
that $C_f=R_f + S_f$, the induced map $R_f \rightarrow f(R_f)$ is finite, and
$f(S_f)$ is a finite set. Let $B_f := f(R_f)_{\rm red}$. By the Zariski-Nagata
theorem on the purity of the branch locus, Theorem X.3.1 [SGA], $B_f$ is a divisor on $X$.

\begin{defin}
A nonsingular model $f: \ Y \rightarrow X$ is said to be monomial if
for every $\beta \in Y$, with $\alpha := f(\beta) \in X$, there exist
regular systems of parameters (r.s.p. for short) $(u,v)$ of
${\cal O}_{X,\alpha}$ and $(x,y)$ of ${\cal O}_{Y,\beta}$ such that
\begin{itemize}
\item[(i)]
If $\alpha \in {\rm Supp}(B_f)$, $B_f$ is locally at $\alpha$
defined by $u=0$ or $uv=0$.
\item[(ii)]
Either
$$
\left\{\matrix{
u & = & \gamma x^a y^b \cr
v & = & \delta x^c y^d \cr } \right. , \leqno (1)
$$
where $\gamma \delta$ is a unit in ${\cal O}_{Y,\beta}$ and $p$ does not
divide $ad-bc$, or
$$
\left\{\matrix{
u & = & \gamma x^a \cr
v & = & \delta x^c \cr } \right. , \leqno (2)
$$
where both of $\gamma \delta$ and
$a \gamma \displaystyle{\partial \delta \over \partial y} -
c \delta {\partial \gamma \over \partial y}$ are units in ${\cal O}_{Y,\beta}$.
\end{itemize}
A monomial model dominating a given model is called a monomial resolution.
\end{defin}

\begin{prop}
A monomial model $f: \ Y \rightarrow X$ has the following properties.
\begin{itemize}
\item[(i)]
$B_f$, $f^*B_f$ and $S_f$ are divisors with strict normal crossings.
\item[(ii)]
For every $\beta \in Y$, with $\alpha := f(\beta) \in X$,
and regular parameters $(u,v)$ in ${\cal O}_{X,\alpha}$ as in {\it (ii)} of
Definition 3,
there exists an affine neighborhood $U$ of $\beta$, and an \'etale
cover $V$ of $U$ such that there are uniformizing parameters
$(\overline x,\overline y)$  on $V$ with
$$
\begin{array}{l}
u=\overline x^a\overline y^b\\
 v =\overline x^c\overline y^d
\end{array}
$$
for some natural numbers $a,b,c,d$ such that $p$ does not divide $ad-bc$.
\end{itemize}
\end{prop}

{\it Proof:} {\it (i)} directly follows from definition 3.

We will prove (ii), under the assumption that case (2) of Definition 3 holds.
There exists an affine neighborhood $U$ of $\beta$
such that $(x,y)$ are uniformizing parameters on $U$, and
$\gamma$, $\delta$, and
$a\gamma\frac{\partial \delta}{\partial y}-c\delta\frac{\partial \gamma}{\partial y}$
are units in $\Gamma(U,{\cal O}_Y)$. $p$ cannot divide both $a$ and $c$.
Without loss of generality, $p$ does not divide $a$.
Set $d=1$, $b=0$,
$$
R=\Gamma(U,{\cal O}_{Y})[\gamma^{\frac{1}{a}},
\delta^{\frac{1}{a}}],
$$
$V= \mbox{Spec}(R)$. $V$ is an \'etale cover of $U$.
Set
$$
\overline
x=x\gamma^{\frac{1}{a}},
\overline y = \delta\gamma^{\frac{-c}{a}}.
$$
$(\overline x,\overline y)$ are uniformizing parameters on $V$ since
$$
\frac{\partial \overline y}{\partial y} =
\frac{1}{a}
\gamma^{\frac{-a}{c}-1}[a\gamma\frac{\partial \delta}{\partial
y}-c\delta \frac{\partial \gamma}{\partial y}]
$$
is a unit in $\Gamma(V,{\cal O}_V)$.\\

Our main result is
\begin{them}
Given a model $f$ of $L/K$, the following properties are equivalent.
\begin{itemize}
\item[(i)]
$f$ admits a minimal (w.r.t. domination) monomial resolution.
\item[(ii)]
$f$ admits a monomial resolution.
\item[(iii)]
$f$ is tamely ramified.
\end{itemize}
\end{them}

Theorem 1 will be proved at the end of section 4. Note that,
if ${\rm char}k=0$, any model of $L/K$ is tamely ramified.\\

{\it Remark:} In case the given model $f: \ Y \rightarrow X$ is {\it finite}, 
there is an easier proof of theorem 1, using Abhyankar's lemma ([Ab5], 2.3.4, [GM]). 
This can be seen as follows; first reduce to $f$ finite and ramified over a 
divisor with strict normal crossings. By Abhyankar's lemma, $f$ can be 
locally described as a Kummer covering, after a local-\'etale change of 
coordinates on $X$. Let $Y'$ be the {\it minimal} resolution of 
singularities of $Y$. 
By explicit computations, it is now seen that $Y' \rightarrow X$ is a monomial 
morphism.\\
From this result, one deduces the {\it existence} of a monomial resolution of 
a given $f$ as in theorem 1, since any model can be dominated by a finite one. 
However, the monomial resolution thus obtained is not in general 
the minimal one.\\

\begin{cor}
Assume that ${\rm char}k=0$ or ${\rm char}k=p>0$ and $p$ does not
divide the degree of the Galois closure $\overline{L} /K$ of $L/K$.
Then any model of $L/K$ admits a minimal monomial resolution.
\end{cor}

{\it Proof:} Assume that ${\rm char}k=p>0$. Let $V$ be a divisorial valuation
ring of $L$ and $\overline{V}$ be a divisorial valuation ring of
$\overline{L}$ such that $V=\overline{V} \cap L$.
Let $W:= V \cap K$ and $e$ and $f$ be the ramification index
and residual degree of the extension $\overline{V}/W$.\\
By V.9.22 of [ZS1], $[\overline{L}:K]=efg$, where $g$ is the number of
conjugates of $\overline{V}$ under the action of ${\rm Gal}(\overline{L} /K)$.
Hence $p$ does not divide $ef$. This implies that $\overline{V}/W$ is
tamely ramified. Consequently, $V/W$ is tamely ramified.\\

\section{The complexity.}

In this section, we only consider nonsingular models $f: \ Y \rightarrow X$
such that both of $B_f$ and $f^*B_f$ are divisors
with strict normal crossings.\\

Let $\alpha$ be a point in $X$. A r.s.p. $(u,v)$ of ${\cal O}_{X,\alpha}$
is said to be {\it admissible} if $\alpha \not \in {\rm Supp}(B_f)$,
or if $\alpha \in {\rm Supp}(B_f)$ and $B_f$ is locally at $\alpha$
defined by $u=0$ or $uv=0$ (see {\it (i)} in definition 3).\\

\begin{defin}
Given a reduced irreducible component $E$ of $S_f$, with
$\alpha:=f(E) \in X$, the {\it complexity} $i_E$ of $E$ is defined
by the following formula:\\
$$
i_E:= \nu_E (S_f) +1 - \max_{(u,v) \ adm.} {\nu_E (uv)} \geq 0,
$$
where $\nu_E$ is the divisorial valuation associated with $E$, and the
maximum is taken over all {\it admissible} r.s.p. at $\alpha$.
\end{defin}

{\it Remark:} it follows from this definition that if $f: \ Y \rightarrow X$
and $f': \  Y' \rightarrow X'$ are two nonsingular models as above, and $E$
(resp. $E'$) is a reduced irreducible component of $S_f$ (resp. $S_{f'}$)
such that
\begin{itemize}
\item[(i)]
${\cal O}_{Y,E} = {\cal O}_{Y',E'}$, and
\item[(ii)]
${\cal O}_{X,\alpha}={\cal O}_{X',{\alpha}'}$, where $\alpha :=f(E)$ and
${\alpha}':=f'(E')$,
\end{itemize}
then $i_E=i_{E'}$. This fact will be repeatedly used in this section.\\

We first recall the following classical birational fact (theorem 3, [Ab1]),
together with its global counterpart (theorem 4.1, [Li]).

\begin{prop}
Let $R$ be a two-dimensional regular local ring with quotient field $K$ and
$S$ be a regular local ring birationally dominating $R$. Assume that $S$ is
either two-dimensional or a divisorial $k$-valuation ring.\\
There exists a unique sequence
$$
R=R_0 \subset R_1 \subset \ldots \subset R_n=S
$$
such that, for $1 \leq i \leq n$, $R_i$ is a quadratic transform of $R_{i-1}$.
\end{prop}

\begin{prop}
Let $R$ be a two-dimensional regular local ring with quotient field $K$ and 
$X \rightarrow {\rm Spec} R$ be a proper birational map with $X$ regular.\\
There exists a sequence
$$
X=X_n \longrightarrow  \ldots \longrightarrow X_1
\longrightarrow X_0={\rm Spec} R
$$
such that, for $1 \leq i \leq n$, $X_i$ is the blow-up of a closed point
of $X_{i-1}$.
\end{prop}

Proposition 2 implies the following.

\begin{cor}
Let $f: \ Y \rightarrow X$ be a nonsingular model as in the beginning of this
section, and let $E$ be a reduced irreducible component of $S_f$,
with $\alpha:=f(E) \in X$. Assume that
$\alpha \not \in {\rm Supp}(B_f)$. Let $\beta \in f^{-1}(\alpha)$. Then,
\begin{itemize}
\item[(i)]
If $\beta \not \in {\rm Supp}(S_f)$, then $M_\alpha {\cal O}_{Y,\beta} =
M_\beta$.
\item[(ii)]
If $\beta \in {\rm Supp}(S_f)$, then $M_\alpha {\cal O}_{Y,\beta}$ is a
principal ideal.
\end{itemize}
\end{cor}

{\it Proof:} Let $\overline{R}$ be the integral closure of
${\cal O}_{X,\alpha}$ in $L$. By X.3.1 [SGA], $\overline{R}$ is a regular
semilocal ring which is unramified over ${\cal O}_{X,\alpha}$. One has 
$\overline{R} \subseteq {\cal O}_{Y,\beta}$ in either case.\\
If $\beta \not \in {\rm Supp}(S_f)$,
${\cal O}_{Y,\beta}=\overline{R}_{M_\beta \cap \overline{R}}$ by proposition
2. Hence $M_\alpha {\cal O}_{Y,\beta} =M_\beta$.\\
If $\beta \in {\rm Supp}(S_f)$, ${\cal O}_{Y,\beta}$ dominates a quadratic
transform $R_1$ of $\overline{R}_{M_\beta \cap \overline{R}}$ by proposition
2. Hence $M_\alpha {\cal O}_{Y,\beta}$ is a principal ideal.\\

\begin{prop}
Let $f: \ Y \rightarrow X$ be a nonsingular model as above,
and let $\beta \in Y$, with $\alpha:=f(\beta)$.
There exists an admissible r.s.p. $(u,v)$ at $\alpha$ such that
for every reduced irreducible component $E$ of $S_f$ passing through $\beta$,
$$
\nu_E (uv)=\max_{(u',v') \ adm.} {\nu_E (u'v')}.
$$
\end{prop}

{\it Proof:}  Since $f^*B_f$ (and hence $S_f$ as well) is a divisor with
strict normal crossings, there exist $s_\beta \leq 2$ components of $S_f$
passing through $\beta$. The above statement is trivial unless $s_\beta =2$,
which we now assume. Since $B_f$ is a divisor with strict normal crossings,
there exist $r_\alpha \leq 2$ components of $B_f$ passing through $\alpha$.
The above statement is trivial if $r_\alpha =2$, and we hence assume
$r_\alpha \leq 1$. Let $E$ be an irreducible component of $S_f$ passing
through $\beta$. We consider two cases:\\

First assume that $r_\alpha =0$. By proposition 2
and Theorem X.3.1 [SGA], there exists a regular local ring $R$, essentially of
finite type and unramified over ${\cal O}_{X,\alpha}$, and a succession of
quadratic transforms
$$
R=R_0 \subset R_1 \subset \ldots \subset R_n ={\cal O}_{Y,\beta},
$$
with $n \geq 1$. Let $u \in {\cal O}_{X,\alpha}$ be a regular parameter.
Let $t_i \in R_i$, $1 \leq i \leq n$, be a regular parameter such that
${\rm ht}\left ( (t_i)\cap R_{i-1}\right )=2$.
Define by induction on $i$, $0 \leq i \leq n$,
elements $u_i \in R_i$ by $u_0=u$, and if $i \geq 1$:
$$
\left\{\matrix{
u_{i-1} & = & t_iu_i & \ {\rm if} \ u_{i-1}R_i \not = R_i \cr
u_{i-1} & = & u_i    & \ {\rm if} \ u_{i-1}R_i      = R_i \cr } \right. .
$$
Let $m_u$, $1 \leq m_u \leq n$, be the largest integer $m$ such that
$u_{m-1}=t_mu_m$. We have:
$$
\nu_E (u)= \sum_{i=1}^{m_u}{\nu_E (t_i)}.
$$
In particular, $\nu_E (u)$ is a non decreasing function of $m_u$. Also
notice that for general $u$, $m_u=1$. A r.s.p. $(u,v)$
satisfying the conclusion of the proposition is then obtained by taking $v$
maximizing $m_v$, and any transversal $u$.\\

Assume now that $r_\alpha =1$. Let $u=0$ be a local equation of $B_f$ at
$\alpha$, and $xy=0$ be a local equation of $(S_f)_{\rm red}$ at $\beta$.
Let $v ,w \in {\cal O}_{X,\alpha}$ be such that
both of $(u,v)$ and $(u,w)$ are (admissible) r.s.p. We have
$$
\left\{\matrix{
u & = & \gamma x^a y^b \cr
v & = & x^c y^d v'\cr
w & = & x^{c'}y^{d'}w' \cr } \right. ,
$$
where $\gamma$ is a unit in ${\cal O}_{Y,\beta}$, $a,b,c,d,c',d'>0$, and
neither $x$ nor $y$ divides $v'w'$.
Assuming that $c'>c$, we will prove that $d' \geq d$ and the conclusion
will follow.\\
 By the Weierstrass preparation theorem, there exists a power
series $P(u) \in \kappa (\alpha)[[u]]$ such that
$$
w = unit \times (v -P(u)) \in {\widehat{{\cal O}}}_{X,\alpha} \simeq
\kappa (\alpha)[[u,v]]. \leqno (3.1)
$$
Let $\lambda \in {\cal O}_{X,\alpha}$ be a unit such that
$$
P(u) \equiv \lambda u^m \ \ {\rm mod} \ u^{m+1},
$$
where $m={\rm ord}_uP$. Since $c'={\rm ord}_xw >c={\rm ord}_xv$, (3.1) implies
that $c=ma$. This gives the congruence
$$
v'y^d \equiv \lambda \gamma^m y^{mb} \ \ {\rm mod} \ x
$$
in ${\cal O}_{Y,\beta}$. Hence $d \leq mb$. It then follows from (3.1) that
$$
d'={\rm ord}_yw \geq {\min} \lbrace {\rm ord}_yv , mb \rbrace =d.
$$
This concludes the proof.\\

Proposition 4 leads to the following definition of the local complexity on $Y$.

\begin{defin}
Let $f: \ Y \rightarrow X$ be a nonsingular model as above, and $\beta \in
{\rm Supp}(S_f)$. The complexity $i_\beta$ of $f$ at $\beta$ is defined by
$$
i_\beta := \max_{E}{i_E} \geq 0,
$$
where the maximum is taken over all reduced irreducible components
of $S_f$ passing through $\beta$.
\end{defin}

\begin{lem}
Let $f: \ Y \rightarrow X$ be a nonsingular model as above, and
$\beta \in {\rm Supp}(R_f)$, with $\alpha:=f(\beta)$. Let $x=0$ be a local
equation of a reduced component $D$ of $R_f$ passing through $\beta$,
and $u=0$ be an equation of $\Delta:=f(D)$ at $\alpha$.
Write $u=x^au'\in {\cal O}_{Y,\beta}$,
where $a \geq 2$, and $x$ does not divide $u'$.

The extension of divisorial valuation rings
${\cal O}_{Y,D} / {\cal O}_{X,\Delta}$ is tamely ramified if and only if
${\rm ord}_DR_f=a-1$.
\end{lem}

{\it Proof:} Choose a r.s.p. $(u,v)$ at $\alpha$, and a r.s.p.
$(x,y)$ at $\beta$. A local equation at $\beta$ of $C_f$ is given by
$$
{\rm Jac}_\beta(f)=x^{a-1} \left ( au'{\partial v \over \partial y} + x
{\rm Jac}(u',v) \right ).
$$
Then ${\rm ord}_DR_f=a-1$ if and only if $p$ does not divide $a$ and
$x$ does not divide $\displaystyle {\partial v \over \partial y}$.

Since ($k$ is perfect) $\alpha$ (resp. $\beta$) is a smooth point of $\Delta$
(resp. $D$) (ex. II.8.1, [Ha]), $\Omega^1_{\Delta /k}$ (resp. $\Omega^1_{D/k}$)
is generated at $\alpha$ (resp. $\beta$) by $dv$ (resp. $dy$).
The inclusion $\kappa({\cal O}_{X,\Delta}) \subseteq \kappa({\cal O}_{Y,D})$
is separable if and only if $\Omega^1_{D / \Delta}$ is a torsion sheaf
(II.8.6.A, [Ha]). Let $\overline{f}: \ D \rightarrow \Delta$ be the finite map
induced by $f$. There is an exact sequence for differentials on $D$
(II.8.11, [Ha])
$$
f^*\Omega^1_{\Delta /k} \buildrel {d\overline{f}} \over {\longrightarrow}
\Omega^1_{D/k} \longrightarrow \Omega^1_{D/\Delta} \longrightarrow 0.
$$
Clearly, $\Omega^1_{D/\Delta}$ is a torsion sheaf if and only if
$(\Omega^1_{D/\Delta})_\beta$ is a torsion ${\cal O}_{Y,\beta}$-module.
The tangent map $d\overline{f}$ is given at $\beta$ by
$$
dv \mapsto \left ( {\partial v \over \partial y} \ {\rm mod} \ x \right ) dy.
$$
It follows that the inclusion $\kappa({\cal O}_{X,\Delta}) \subseteq
\kappa({\cal O}_{Y,D})$ is separable if and only if $x$ does not divide
$\displaystyle {\partial v \over \partial y}$.

Summing up, ${\rm ord}_DR_f=a-1$ if and only the extension of divisorial
valuation rings ${\cal O}_{Y,D} / {\cal O}_{X,\Delta}$
is tamely ramified as required.\\

The following proposition characterizes a monomial model by way of its
maximal complexity.\\

\begin{prop}
Let $f: \ Y \rightarrow X$ be a nonsingular model as above, and let
$\beta \in Y$. $f$ is monomial at $\beta$ (i.e. has the local form (1)
or (2) of definition 3) w.r.t. some admissible r.s.p. at
$\alpha :=f(\beta)$ if and only if
\begin{itemize}
\item[(i)]
For every irreducible component $D$ of $R_f$ passing through $\beta$,
with $\Delta :=f(D)$, the extension of divisorial valuation rings
${\cal O}_{Y,D} / {\cal O}_{X,\Delta}$ is tamely ramified.
\item[(ii)]
$i_\beta = 0$ if $\beta \in {\rm Supp}(S_f)$.
\end{itemize}
In particular, $f$ is a monomial model if and only if {\it (i)} holds
for all components of $R_f$ and {\it (ii)} holds for all 
$\beta \in {\rm Supp}(S_f)$.
\end{prop}

{\it Proof:} Choose an admissible r.s.p. $(u,v)$ at $\alpha$. If 
$\beta \in {\rm Supp}(S_f)$, assume furthermore that $(u,v)$ achieves 
$i_\beta=0$ (proposition 4).
We first prove the if part. We consider six cases.\\

{\it Case 1.} $\beta \not \in {\rm Supp}(C_f)$. Consequently
$M_\alpha {\cal O}_{Y,\beta} = M_\beta$.\\

For cases 2 to 6, assume in addition that $(x,y)$ is chosen such that
$(C_f)_{\rm red}$ has local equation $x=0$ or $xy=0$ at $\beta$.\\

{\it Case 2.} $\beta$ {\it is a smooth point of} ${\rm Supp}(R_f)$ {\it and}
$\beta \not \in {\rm Supp}(S_f)$. Then $x=0$ is a local equation
of $D:=(R_f)_{\rm red}$ at $\beta$ and, say, $u=0$ is a local equation
of $f(D)$ at $\alpha$. We have $u=x^au'$, where $a \geq 2$ and $x$
does not divide $u'$. The local equation at $\beta$ of $C_f$ is given by
$$
{\rm Jac}_\beta(f)=x^{a-1} \left ( au'{\partial v \over \partial y} + x
{\rm Jac}(u',v) \right ). \leqno (3.2)
$$
By lemma 1, ${\rm ord}_DR_f = a-1$. Since $C_f$ is a divisor with strict
normal crossings, $au'{\partial v \over \partial y}$ is a unit. $f$ is then
reduced at $\beta$ to the monomial form (1)
$$
\left\{\matrix{
u & = & \gamma x^a  \cr
v & = & y \cr } \right. ,
$$
where $\gamma$ is a unit and $p$ does not divide $a$.\\

{\it Case 3.} $\beta$ {\it is a singular point of} ${\rm Supp}(R_f)$.
Then $xy=0$ is a local equation of $D_1+D_2=(R_f)_{\rm red}$ at $\beta$.
Suppose, if possible, that $u=0$ is a local equation of
$f(D_1 \cup D_2)$ at $\alpha$. Hence
$u=x^ay^bu'$, where $a,b \geq 2$ and neither $x$ nor $y$ divides $u'$.
The local equation at $\beta$ of $C_f$ is given by
$$
{\rm Jac}_\beta(f)=x^{a-1} y^{b-1} \left ( ayu'{\partial v \over \partial y}
-bxu'{\partial v \over \partial x} + xy {\rm Jac}(u',v) \right ).
$$
By lemma 1, ${\rm ord}_{D_1}R_f = a-1$ and ${\rm ord}_{D_2}R_f = b-1$.
Since $C_f$ is a divisor with strict normal crossings, one gets that
$$
ayu'{\partial v \over \partial y}-bxu'{\partial v \over \partial x} +
xy {\rm Jac}(u',v)
$$
is a unit: a contradiction. So $uv=0$ is a local equation of
$f(D_1 \cup D_2)$ at $\alpha$, and $f$ is then reduced at $\beta$ to the
monomial form (1)
$$
\left\{\matrix{
u & = & \gamma x^a  \cr
v & = & \delta y^d \cr } \right. .
$$
The local equation of $C_f$ at $\beta$ is given by
$$
{\rm Jac}_\beta(f)=x^{a-1} y^{d-1} (ad \gamma \delta +g), \leqno (3.3)
$$
where $g$ is a nonunit. Hence $\gamma \delta$ is a unit and $p$ does not
divide $ad$.\\

{\it Case 4.} $\beta$ {\it is a smooth point of} ${\rm Supp}(S_f)$
{\it and} $\beta \not \in {\rm Supp}(R_f)$. Then $x=0$ is a local equation
of $E:=(S_f)_{\rm red}$ at $\beta$. We have $u=x^au'$ and $v=x^cv'$, where
$a,c \geq 1$ and $x$ does not divide $u'v'$. The local equation at $\beta$
of $C_f$ is given by
$$
{\rm Jac}_\beta(f)=x^{a+c-1} \left ( au'{\partial v' \over \partial y}
-cv'{\partial u' \over \partial y} + x {\rm Jac}(u',v') \right ). \leqno (3.4)
$$
By assumption {\it (ii)}, ${\rm ord}_ES_f=a+c-1$. Since $C_f$ is
a divisor with strict normal crossings, it follows that
$au'{\partial v' \over \partial y}-cv'{\partial u' \over \partial y}$
is a unit.\\
Suppose $u'v'$ is not a unit, say, $v'$ is not. Hence
$au'{\partial v' \over \partial y}$ is a unit. $f$ is then reduced at
$\beta$ to the monomial form (1)
$$
\left\{\matrix{
u & = & \gamma x^a  \cr
v & = & x^cy \cr } \right. ,
$$
where $\gamma$ is a unit and $p$ does not divide $a$.\\
Suppose $u'v'$ is a unit. $f$ is then reduced at $\beta$ to the
monomial form (2).\\

{\it Case 5.} $\beta \in {\rm Supp}(R_f) \cap {\rm Supp}(S_f)$. Then $x=0$
(resp. $y=0$) is a local equation of $D:=(R_f)_{\rm red}$ (resp.
$E:=(S_f)_{\rm red}$) at $\beta$. We have $u=x^ay^bu'$ and $v=y^dv'$, where
$a \geq 2$, $b,d \geq 1$ and neither $x$ nor $y$ divides $u'v'$.
The local equation at $\beta$ of $C_f$ is given by
$$
{\rm Jac}_\beta(f)=x^{a-1} y^{b+d-1} \left ( adu'v'
+ayu'{\partial v' \over \partial y} + x g \right )   \leqno (3.5)
$$
for some $g$. By lemma 1, ${\rm ord}_DR_f = a-1$. By assumption {\it (ii)},
${\rm ord}_ES_f = b+d-1$. Since $C_f$ is a divisor with strict normal
crossings, one gets that $adu'v'$ is a unit.
$f$ is then reduced at $\beta$ to the monomial form (1)
$$
\left\{\matrix{
u & = & \gamma x^ay^b  \cr
v & = & \delta y^d \cr } \right. ,
$$
where $\gamma \delta$ is a unit and $p$ does not divide $ad$.\\

{\it Case 6.} $\beta$ {\it is a singular point of} ${\rm Supp}(S_f)$.
Then $xy=0$ is a local equation of $E_1+E_2=(S_f)_{\rm red}$ at $\beta$.
Hence $u=x^ay^bu'$ and $v=x^cy^dv'$, where $a,b,c,d \geq 1$ and neither $x$
nor $y$ divides $u'v'$. The local equation at $\beta$ of $C_f$ is given by
$$
{\rm Jac}_\beta(f)=x^{a+c-1} y^{b+d-1} \left ( (ad-bc)u'v'+g \right ),
\leqno (3.6)
$$
where $g$ is a nonunit. By assumption {\it (ii)}, ${\rm ord}_{E_1}S_f = a+c-1$
and ${\rm ord}_{E_2}S_f = b+d-1$. Since $C_f$ is a divisor with strict normal
crossings, one gets that $(ad-bc)u'v'$ is a unit. $f$ is then reduced at
$\beta$ to the monomial form (1).\\

The only if part of the proposition easily follows by applying formulas
(3.2) to (3.6) to the monomial expression (1) or (2). The last statement
is obvious. This completes the proof.\\

In order to construct a monomial model dominating a given model as above,
it is necessary to study the behaviour of the complexity $i_\beta$ under
blow-up. This is achieved in propositions 6 and 7 below.

\begin{prop}
Let $f: \ Y \rightarrow X$ be a nonsingular model as above, and let $E$ be
a reduced irreducible component of $S_f$, with $\alpha :=f(E) \in X$.\\
Assume that $M_\alpha {\cal O}_Y$ is locally principal. Let
$\eta : \ X' \rightarrow X$ be the blowing-up of $\alpha$, and
$f': \ Y \rightarrow X'$ be the induced map.\\
Assume in addition that ${\alpha}':=f'(E) \in X'$ is a point,
and let $i_E$ (resp. $i'_E$) be the complexity of $E$ w.r.t. $f$
(resp.  $f'$). Then $i'_E \leq i_E$.
\end{prop}

{\it Proof:} Pick an admissible r.s.p. $(u,v)$ at $\alpha$ achieving $i_E$.
Say, $\nu_E(u) = \min_{t \in M_\alpha} \lbrace \nu_E (t) \rbrace$. Then 
$u=0$ is a local equation of the exceptional divisor of $\eta$ at ${\alpha}'$. 
Pick $v'$ such that $(u,v')$ is an admissible r.s.p. at ${\alpha}'$ with 
$\nu_E (v')$ maximal. Consequently
$$
\left.\matrix{
i'_E & \leq & \nu_E(S_{f'})        & + & 1 & - & \nu_E(uv') & \cr
     &  =   & \nu_E(S_f) -\nu_E(u) & + & 1 & - & \nu_E(uv') &
= i_E +\nu_E(v) -\nu_E(uv') \cr} \right. .
$$
If $\nu_E(u)=\nu_E(v)$, then $\nu_E(v) -\nu_E(uv') < 0$ and $i'_E < i_E$.\\
If $\nu_E(u)<\nu_E(v)$, then $(u,{v \over u})$ is an admissible r.s.p. at
${\alpha}'$ and consequently $\nu_E(v)=\nu_E(u{v \over u}) \leq \nu_E(uv')$,
i.e. $i'_E \leq i_E$.\\

\begin{lem}
Let $f: \ Y \rightarrow X$ be a nonsingular model as above, and
$\beta \in {\rm Supp}(S_f)$, with $\alpha:=f(\beta)$. Let $x=0$ be a local
equation of a reduced component $E$ of $S_f$ passing through $\beta$,
and $(u,v)$ be a r.s.p. at $\alpha$ achieving $i_E$. Write
$$
\left\{\matrix{
u & = & x^au'  \cr
v & = & x^cv' \cr } \right. ,
$$
where $a,c \geq 1$, and $x$ does not divide $u'v'$.
Let $\delta := {\rm g.c.d.}(a,c)$.\\
Assume that ${u'}^{c \over \delta}$ divides ${v'}^{a \over \delta}$, and
that ${\displaystyle{{v'}^{a \over \delta} \over {u'}^{c \over \delta}}}$ is
not a unit in ${\cal O}_{Y,\beta}$.\\
The extension of divisorial valuation rings
${\cal O}_{Y,E} / {\cal O}_{Y,E}\cap K$ is tamely ramified
if and only if $i_E=0$.
\end{lem}

{\it Proof:} Choose a r.s.p. $(x,y)$ at $\beta$.
A local equation at $\beta$ of $C_f$ is given as in (3.4) by
$$
{\rm Jac}_\beta(f)=x^{a+c-1} \left ( au'{\partial v' \over \partial y} -
cv'{\partial u' \over \partial y} + x{\rm Jac}(u',v') \right ). \leqno (3.7)
$$
Then $i_E=0$ if and only if $x$ does not divide
$au'{\partial v' \over \partial y} - cv'{\partial u' \over \partial y}$.
Let $\varphi:= {\displaystyle {v^{a \over \delta} \over u^{c \over \delta}}}$.
By assumption, $\varphi \in {\cal O}_{Y,\beta}$ and $\varphi$ is not a unit.\\

Let $I$ be the integral closure of the ideal
$(u^{c \over \delta}, v^{a \over \delta})$. Then $I$ is a {\it simple}
complete $M_\alpha$-primary ideal (p.385, appendix 5 of [ZS2]). 
There are local inclusions
$$
{\cal O}_{X,\alpha} \subset R_{\cal Q} \subseteq {\cal O}_{Y,\beta} \ ,
\leqno (3.8)
$$
where $R=\displaystyle{{\cal O}_{X,\alpha}
\left [{I \over u^{c \over \delta}} \right ]}$
and ${\cal Q}:=M_\beta \cap R$. By construction, $\varphi \in R$ and is
neither a unit nor is divisible by $x$ in ${\cal O}_{Y,\beta}$.
This implies that ${\rm ht}((x) \cap R)=1$.\\

{\it Remark:} In case $E$ is the unique component of $S_f$ passing through 
$\beta$, the ring $R_{\cal Q}$ is the local ring {\it lying below} 
${\cal O}_{Y,\beta}$ according to Abhyankar's terminology (cf. prop. 2 and 
def. 4 of [Ab4]).\\

By Zariski's theory of complete ideals in two-dimensional regular local
rings ((E) p.391, [ZS2]), there is a 1-1 correspondance between simple complete
$M_\alpha$-primary ideals of ${\cal O}_{X,\alpha}$ and divisorial valuation
rings of $K$ dominating ${\cal O}_{X,\alpha}$;
the reduced exceptional divisor of the
blow-up $\overline{X}:={\rm Proj}({\bigoplus}_{n \geq 0}I^n) \rightarrow
{\rm Spec}{\cal O}_{X,\alpha}$ is an irreducible curve $F$ and
$V:={\cal O}_{\overline{X},F}$ is a divisorial valuation ring (proposition
21.3 and remark following, [Li]). By what preceeds,
$$
V = R_{(x) \cap R} = {\cal O}_{Y,E} \cap K.
$$

Let $t$ be a uniformizing parameter of $V$. Since $I$ is a monomial ideal,
the value group of $V$ is generated by the values ${\rm ord}_tu$ and
${\rm ord}_tv$; this follows from [Sp], lemmas 8.1 and 8.2, and corollary 8.5,
where $k$ needs not be algebraically closed in the special case of a
monomial ideal.  Since ${\rm ord}_t\varphi = 0$, this implies that
${\rm ord}_tu={a \over \delta}$ and ${\rm ord}_tv = {c \over \delta}$.
Hence $IV=(t)^{{a \over \delta}{c \over \delta}}$. On the other hand,
$I{\cal O}_{Y,E}= (x)^{{ac \over \delta}}$. Hence the ramification index of
the extension of divisorial valuation rings ${\cal O}_{Y,E} / V$
is equal to $\delta$.\\

The ideal $I$ is generated by all monomials $u^mv^n$ such that
$$
{m \over {c \over \delta}} + {n \over {a \over \delta}} \geq 1.
$$
Since ${\rm g.c.d.}({a \over \delta},{c \over \delta})=1$, one gets that
$$
\nu_E(u^mv^n)= ma+nc > {ac \over \delta}= \nu_E(u^{c \over \delta})
$$
for all such monomials provided $(m,n) \not \in \lbrace ({c \over \delta},0),
(0,{a \over \delta}) \rbrace$. This proves that
$$
{R \over (x) \cap R}=\kappa(\alpha)[\overline{\varphi}],
$$
where $\overline{\varphi}$ is the image of $\varphi$ in the ring to the
left. Let $\overline{\alpha}$ be the point of $F$ corresponding to ${\cal Q}$.
Then $\overline{\varphi}$ is a regular parameter of
${\cal O}_{F,\overline{\alpha}}$. By (3.8), the rational map
$Y \cdots \rightarrow \overline{X}$ is defined at $\beta$. Besides,
$\beta$ (resp. $\overline{\alpha}$) is a smooth point of $E$ (resp. $F$).
Arguing as in lemma 1, one deduces that
the residue field extension $\kappa({\cal O}_{Y,E})/\kappa(V) \ (=K(E)/K(F))$
is separable if and only if $x$ does not divide
${\partial \varphi \over \partial y}$ in ${\cal O}_{Y,\beta}$. We have
$$
u'v'{\partial \varphi \over \partial y}=  \varphi
\left ( {a \over \delta}u' {\partial v' \over \partial y} -
{c \over \delta} v'{\partial u' \over \partial y} \right ).
$$

Summing up, ${\cal O}_{Y,E} / V$ is tamely ramified if and only if
$x$ does not divide $au'{\partial v' \over \partial y} -
cv'{\partial u' \over \partial y}$.
By (3.7), this is equivalent to $i_E=0$.\\

\begin{prop}
Let $f: \ Y \rightarrow X$ be a nonsingular model as above, and let
$\beta \in Y$, with $\alpha:=f(\beta) \in X$. Assume that
$M_\alpha {\cal O}_{Y,\beta}$ is not a principal ideal.
Let $Y' \rightarrow Y$ be the blowing-up of $\beta$, with
exceptional divisor $E'$, and $f': \ Y' \rightarrow X$ be the composed map.
Let $i_{E'}$ be the complexity of $E'$ w.r.t. $f'$.\\
Assume in addition that either $f$ is monomial at $\beta$ (i.e. has the
local form (1) or (2) of definition 3) w.r.t. some admissible r.s.p. at
$\alpha$, or that the map $f^{(\alpha)}: \ Y^{(\alpha)} \rightarrow
{\rm Spec}{\cal O}_{X,\alpha}$ obtained from
$f$ by the base change ${\rm Spec}{\cal O}_{X,\alpha} \rightarrow X$
is tamely ramified (definition 2). The following holds
\begin{itemize}
\item[(i)]
if $\beta \not \in {\rm Supp}(S_f)$, then $i_{E'}=0$.
\item[(ii)]
if $\beta \in {\rm Supp}(S_f)$, then $i_{E'} \leq i_\beta$,
and $i_{E'} < i_\beta$ if $i_\beta >0$.
\end{itemize}
\end{prop}

{\it Proof:} First assume that $\alpha \not \in {\rm Supp}(B_f)$.
By corollary 2,
this implies that $M_\alpha {\cal O}_{Y,\beta} = M_\beta$, since
$M_\alpha {\cal O}_{Y,\beta}$ is not a principal ideal. We have
$i_{E'}=1+1-2=0$ in this case.\\

Now assume that $\alpha \in {\rm Supp}(B_f)$. Pick an admissible r.s.p.
$(u,v)$ at $\alpha$ achieving $i_E$ for every reduced irreducible component
$E$ of $S_f$ passing through $\beta$. We consider six cases as in
proposition 5.\\

{\it Case 1.} We have $M_\alpha {\cal O}_{Y,\beta} = M_\beta$
and $i_{E'}=0$ as above.\\

{\it Case 2.} By definition or by proposition 5, $f$ is reduced at $\beta$
to the monomial form (1)
$$
\left\{\matrix{
u & = & \gamma x^a  \cr
v & = & y \cr } \right. ,
$$
where $\gamma$ is a unit and $p$ does not divide $a$.
Hence $i_{E'}=a+1-(a+1)=0$.\\

{\it Case 3.} By definition or by proposition 5, $f$ is reduced at $\beta$
to the monomial form (1)
$$
\left\{\matrix{
u & = & \gamma x^a  \cr
v & = & \delta y^d \cr } \right. ,
$$
where $\gamma \delta$ is a unit and $p$ does not divide $ad$.
Hence $i_{E'}= a+d-1+1-(a+d)=0$.\\

{\it Case 4.} First assume that $f$ is monomial at $\beta$. By proposition 5,
$i_\beta=0$ and $f$ is reduced at $\beta$ to the monomial form (1)
$$
\left\{\matrix{
u & = & \gamma x^a  \cr
v & = & x^cy \cr } \right. ,
$$
where $\gamma$ is a unit and $p$ does not divide $a$.
Hence $i_{E'}=a+c +1 -(a+c+1)=0$.\\
Assume now that $f^{(\alpha)}$ is tamely ramified. Write
$$
\left\{\matrix{
u & = & x^au'  \cr
v & = & x^cv' \cr } \right. ,
$$
where $a,c \geq 1$ and $x$ does not divide $u'v'$. Then
$i_{E'} \leq i_\beta +1 -{\rm ord}_{\beta}(u'v')$.
Since $M_\alpha {\cal O}_{Y,\beta}$ is not a principal ideal,
$u'v'$ is not a unit.
Hence $i_{E'} \leq i_\beta$. Assume equality holds. Then $u'$ or $v'$ is a
unit, say $u'$, and $v'$ is not. Hence lemma 2 applies, and gives
that $i_\beta =0$.\\

{\it Case 5.} First assume that $f$ is monomial at $\beta$. By proposition 5,
$i_\beta=0$ and $f$ is reduced at $\beta$ to the monomial form (1)
$$
\left\{\matrix{
u & = & \gamma x^ay^b  \cr
v & = & \delta y^d \cr } \right. ,
$$
where $\gamma \delta$ is a unit and $p$ does not divide $ad$.
Hence $i_{E'}=a+b+d-1 +1 -(a+b+d)=0$.\\
Assume now that $f^{(\alpha)}$ is tamely ramified. Write
$$
\left\{\matrix{
u & = & x^au'  \cr
v & = & x^cy^dv' \cr } \right. ,
$$
where $a,c \geq 1$, $d \geq 2$, and neither $x$ nor $y$ divides $u'v'$. Let
$D$ be the reduced component of $R_f$ with equation $y=0$. By
lemma 1, ${\rm ord}_DR_f =d-1$.
Then $i_{E'} \leq i_\beta -{\rm ord}_{\beta}(u'v') \leq i_\beta$.
Assume $i_{E'} = i_\beta$. Then $u'v'$ is a unit.
By lemma 2, this implies that $i_\beta =0$.\\

{\it Case 6.} First assume that $f$ is monomial at $\beta$. By proposition 5,
$i_\beta=0$ and $f$ is reduced at $\beta$ to the monomial form (1)
$$
\left\{\matrix{
u & = & \gamma x^ay^b  \cr
v & = & \delta x^cy^d \cr } \right. ,
$$
where $\gamma \delta$ is a unit and $p$ does not divide $ad-bc$. Hence
$i_{E'}=a+b+c+d-1 +1 -(a+b+c+d)=0$.\\
Assume now that $f^{(\alpha)}$ is tamely ramified. Write
$$
\left\{\matrix{
u & = & x^ay^bu'  \cr
v & = & x^cy^dv' \cr } \right. ,
$$
where $a,b,c,d \geq 1$, and neither $x$ nor $y$ divides $u'v'$. Let $E_1$
(resp. $E_2$) be the reduced component of $S_f$ with equation $x=0$
(resp. $y=0$). Then
$$
i_{E'} \leq i_{E_1}+i_{E_2} -{\rm ord}_{\beta}(u'v'). \leqno (3.9)
$$
Since $\alpha \in {\rm Supp}(B_f)$ and $(u,v)$ is admissible, $u=0$ is a
local equation of a component of ${\rm Supp}(B_f)$. Since $f^*B_f$ is a
divisor with strict normal crossings, $u'$ is a unit. \\
Suppose that $v'$ is not a unit. By possibly permuting $x$ and $y$, it can be
assumed that $ad-bc \geq 0$. Lemma 2 hence applies, and we get $i_{E_1}=0$.
Hence $i_\beta = i_{E_2}$. By (3.9), this gives $i_{E'} < i_\beta$.\\
Suppose that $v'$ is a unit. Since $M_\alpha {\cal O}_{Y,\beta}$ is not a
principal ideal, $ad-bc \not =0$.
After possibly permuting $u$ and $v$, and $x$ and $y$,
lemma 2 applies w.r.t. both of $E_1$ and $E_2$.
Hence $i_\beta = i_{E_1}=i_{E_2}=0$ and this gives $i_{E'}=0$ by (3.9).\\

\section{The algorithm}

Let $f: Y \rightarrow X$ be a {\it nonsingular} model of $L/K$, and $\alpha$ a
point in $X$. We define a new nonsingular model $f_\alpha$ dominating $f$
as follows: let $X_\alpha \rightarrow X$ be the blowing-up
of $\alpha$, and $Y_\alpha \rightarrow Y$ be the minimal
composition of point blowing-ups such that $M_\alpha {\cal O}_{Y_\alpha}$
is locally invertible (i.e. the minimal resolution of singularities of
the blow-up of $Y$ along the ideal $M_\alpha {\cal O}_Y$).
By the universal property of blow-up,
there exists a map $f_\alpha : \ Y_\alpha \rightarrow X_\alpha$.

\begin{lem}
The above model $f_\alpha$ is the minimal (w.r.t. domination) nonsingular
model $f': Y' \rightarrow X'$ of $L/K$ dominating $f$, and such that the
center on $X'$ of the $M_\alpha$-adic valuation of $K$ is a curve.
\end{lem}

{\it Proof:} Let $f': \ Y' \rightarrow X'$ be a nonsingular model dominating
$f$ such that the center on $X'$ of the $M_\alpha$-adic valuation of $K$
is a curve. By proposition 3, $X' \rightarrow X$ factors through the
blow-up $X_\alpha$ of $X$ at $\alpha$. Since $M_\alpha {\cal O}_{Y'}$ is
locally invertible, $Y' \rightarrow Y$ factors through $Y_\alpha$.\\

\begin{prop}
There exists a minimal nonsingular model $\buildrel \sim \over {f}$
dominating any given model $f$ of $L/K$, and such that both of
$B_{\buildrel \sim \over {f}}$ and
${\buildrel \sim \over {f}}^* B_{\buildrel \sim \over {f}}$
are divisors with strict normal crossings. Any nonsingular model
dominating $\buildrel \sim \over {f}$ has the same property.
\end{prop}

{\it Proof:} Since any algebraic $k$-surface admits a minimal resolution of
singularities (A. p.155, [Li2]), $X$ can be replaced with its minimal 
resolution $X'$ and
$Y$ by the minimal resolution $Y'$ of the normalization of $X$ in $L$.
It hence can be assumed that $f: \ Y \rightarrow X$
is nonsingular. Let $f': \ Y' \rightarrow X'$ be a nonsingular model
dominating $f$. Let $\alpha \in X$ be a point of ${\rm Supp}(B_f)$ which
is not a strict normal crossing. If $B_{f'}$ is a divisor with strict
normal crossings, $f'$ is not an isomorphism above $\alpha$. By lemma 3 and
proposition 3, $f'$ dominates $f_\alpha$.\\
By embedded resolution of curves in surfaces (V.3.9, [Ha]), 
we hence may assume that
$B_f$ has strict normal crossings. Suppose that $f^*B_f$ does not have only
strict normal crossings. Let $\eta: \ \buildrel \sim \over {Y} \rightarrow Y$
be the minimal composition of point blow-ups such that $\eta^*C_f$ has only
strict normal crossings.
Set $\buildrel \sim \over {f}: \ \buildrel \sim \over {Y} \rightarrow X$
to be the morphism induced by $f$.\\

\begin{lem}
A monomial model $f: Y \rightarrow X$ is tamely ramified.
\end{lem}

{\it Proof:} Let $V$ be a divisorial $k$-valuation ring of $L$ having
a center in $Y$.\\
{\it Case 1: the center of $V$ in $X$ is a curve}. By proposition 5{\it (i)},
$V / V \cap K$ is tamely ramified.\\
{\it Case 2 : the center of $V$ in $X$ is a point $\alpha$}. Consider the
model $f_\alpha: \ Y_\alpha \rightarrow X_\alpha$. By propositions 5, 6 and 7,
$f_\alpha$ also is a monomial model. If the center of $V$ in $X_\alpha$ is
a curve, we are done by case 1. If it is a point $\alpha_1$, then
${\cal O}_{X, \alpha} \subset {\cal O}_{X_\alpha, \alpha_1}$, and we apply
again case 2 to $f_{\alpha_1}$.\\
This process terminates after a finite number of steps by proposition 2
applied to $R={\cal O}_{X, \alpha}$ and $S=V \cap K$.\\

{\it Proof of theorem 1 (stated at the end of section 2):}
$(i) \Longrightarrow (ii)$ is trivial and $(ii) \Longrightarrow (iii)$ has
been proved in lemma 4 above. We prove $(iii) \Longrightarrow (i)$.\\

By proposition 8, it can be assumed that $f$ is a nonsingular model and that
both of $B_f$ and $f^*B_f$ are divisors with strict normal crossings. Hence
the results of section 3 apply. \\

{\bf The algorithm:} Assume furthermore that $f$ is {\it not} a monomial model.
By proposition 5, there exists a reduced irreducible component $E$ of
$S_f$ with $i_E>0$. Choose such an $E$ with $i_E$ maximal, and let
$\alpha :=f(E) \in X$. We get a new nonsingular model
$f_\alpha: \  Y_\alpha \rightarrow X_\alpha$ dominating $f$ and such that
both of $B_{f_\alpha}$ and $f_\alpha^*B_{f_\alpha}$ are divisors with
strict normal crossings by proposition 8 above.\\
Iterate the process if $f_\alpha$ is not monomial. This gives rise to a
sequence of nonsingular models $f,f_{\alpha_1}, \ldots , f_{\alpha_i}, \ldots$,
such that $f_{\alpha_i}$ dominates $f_{\alpha_{i-1}}$ for $i \geq 1$. It
will be proved below that for some $i \geq 1$, $f_{\alpha_i}$ is the
minimal monomial resolution of $f$.\\

{\it Proof of theorem 1 continued:} Let $E$ be a reduced irreducible
component of $S_f$ such that $i_E>0$. Let $\alpha := f(E)$.
We first claim that any monomial model $f': \ Y' \rightarrow X'$
dominating $f$ (if there exists one) dominates $f_\alpha$ as well.
By lemma 3, it is sufficient to show that the
center on $X'$ of the $M_\alpha$-adic valuation of $K$ is a curve. Since $X'$
is nonsingular, it is also sufficient by proposition 3 to prove that
$X' \rightarrow X$ is not an isomorphism above $\alpha$.
Assume the contrary. Let $f^{(\alpha)}$
(resp. ${f'}^{(\alpha)}$) be the map obtained from $f$ (resp. $f'$) by the
base change ${\rm Spec}{\cal O}_{X,\alpha} \hookrightarrow X$. There is a
commutative diagram with proper maps
$$
\diagram{
{Y'}^{(\alpha)} & \hfr{{f'}^{(\alpha)}}{}& {\rm Spec}{\cal O}_{X,\alpha}\cr
\vfd{\pi}{} & & \vfd{}{\kern-2pt\displaystyle\wr} \cr
Y^{(\alpha)} & \hfr{f^{(\alpha)}}{}& {\rm Spec}{\cal O}_{X,\alpha} \cr} .
$$
Let $E'$ be the strict transform of $E$ in ${Y'}^{(\alpha)}$. By definition,
$i_{E'} = i_E$. But ${f'}^{(\alpha)}$ is monomial by assumption and
thus $i_{E'}=0$ by proposition 5. This is a contradiction, since $i_E >0$,
and the claim is proved.\\

Let
$$
I_f:=\max_{\beta \in {\rm Supp}(S_f)}i_\beta \ ,
$$
and
$$
\Sigma_f:= \lbrace {\cal O}_{Y,E} \mid E \ {\rm is} \ {\rm a} \ {\rm reduced}
\ {\rm irreducible} \ {\rm component} \ {\rm of} \ S_f \ {\rm with} \
i_E=I_f \rbrace.
$$

To conclude the proof, it must be shown that for some $i \geq 1$, the
model $f_{\alpha_i}$ in the algorithm above is monomial, i.e.
$I_{f_{\alpha_i}}=0$. Assume not. By lemma 5 below, $(I_{f_{\alpha_j}},
\Sigma_{f_{\alpha_j}})$ is constant for large enough $j$.
Pick a divisorial valuation ring
$V \in \displaystyle{\bigcap_{j \geq 0}\Sigma_{f_{\alpha_j}}}$ of $L$, such
that for infinitely many values of $j$, $\alpha_j$ is the center of $V$ in
$X_{\alpha_j}$. This gives rise to an increasing sequence of 
quadratic transforms
$({\cal O}_{X_{\alpha_j},\alpha_j})$ dominated by $V \cap K$. But any such
sequence must be finite by proposition 2.\\

\begin{lem}
With notations as above, assume $I_f >0$ and let $E$ be a reduced irreducible
component of $S_f$ with $i_E=I_f$. Let $\alpha := f(E)$.

Then $(I_{f_\alpha},\Sigma_{f_{\alpha}}) \leq (I_f,\Sigma_f)$
for the lexicographical ordering, where the second summand is (partially)
ordered by inclusion.
\end{lem}

{\it Proof:} Let $f': \ Y_\alpha \rightarrow X$. By proposition 7,
$(I_{f'},\Sigma_{f'}) = (I_f,\Sigma_f)$. By proposition 6,
$(I_{f_\alpha},\Sigma_{f_{\alpha}}) \leq (I_{f'},\Sigma_{f'})$.\\

\section{Toroidalization of morphisms  of surfaces}

Suppose that $k$ is an algebraically closed field of characteristic zero. We
recall the definitions of
toroidal varieties and morphisms from [KKMS] and [AKa].

Suppose that $X$ is a normal $k$-variety, with  an open subset $U_X$. The
embedding $U_X\subset X$
is {\it toroidal} if for every $x\in X$ there exists an affine  toric variety
$X_{\sigma}$, a point $s\in X_{\sigma}$,
and an isomorphism  $\hat{\cal O}_{X,x}\cong \hat{\cal O}_{X_{\sigma,s}}$
such that the ideal of $X-U_X$ corresponds to the ideal of $X_{\sigma}-T$,
where $T$ is the torus in $X_{\sigma}$.
 Such a pair
$(X_{\sigma},s)$ is called a local model at $x\in X$.

A dominant morphism $f:(U_X\subset X)\rightarrow (U_B\subset B)$ of toroidal
embeddings is called {\it toroidal} if for every closed point $x\in X$ there
exist local models $(X_{\sigma},s)$ at $x$,
$(X_{\tau},t)$ at $f(x)$ and a toric morphism $g:X_{\sigma}\rightarrow
X_{\tau}$ such that the following diagram commutes
$$
\begin{array}{lll}
\hat{\cal O}_{X,x} & \leftarrow & \hat{\cal O}_{X_{\sigma},s}\\
\hat f^*\uparrow&&\uparrow \hat g^*\\
\hat{\cal O}_{B,f(x)} & \leftarrow & \hat{\cal O}_{X_{\tau},t}
\end{array}
$$

By a $k$-surface, we mean a proper, 2 dimensional, integral, normal
$k$-variety.

Suppose that $X$ is a nonsingular $k$-surface, and $D_X$ is a SNC (Simple
Normal Crossings)
 divisor on $X$. Then the embedding
$X-D_X\subset X$ is toroidal.

In this section, we will consider morphisms $f:Y\rightarrow X$, where $X$ and
$Y$ are
 $k$-surfaces with respective (Weil) divisors
$D_X$ and $D_Y$ such that $f^{-1}(D_X)=D_Y$, set theoretically. If $\eta :
X_1\rightarrow X$ is a birational
proper morphism of  $k$-surfaces, we can define a  divisor
$D_{X_1}=\eta^{-1}(D_X)$ on $X_1$.
If $D_X$ is a SNC divisor, and $X_1$ is nonsingular, then $D_{X_1}$ is a SNC
divisor.

We will say that $f:Y\rightarrow X$ is toroidal relative to $D_Y$ and $D_X$ if
$$
f:(Y-D_Y\subset Y)\rightarrow (X-D_X\subset X)
$$
 is toroidal.

We will  prove that morphisms of $k$-surfaces can be made toroidal. While this
result is known to be true,
for instance it is implicit
in [AKi], the result is of sufficient interest that we give a statement of the
theorem, and an outline of a proof.
 Recall that, in this section, $k$ is an
algebraically closed field of characteristic zero.

\begin{them} Suppose that $f:Y\rightarrow X$ is a dominant morphism of
$k$-surfaces, and $D_X$,
 $D_Y$ are respective divisors on $X$ and $Y$, such that $f^{-1}(D_X)=D_Y$ and
$D_Y$ contains all
 singular points of $f$ and of $Y$. Then
there exist  projective birational morphisms
$\beta:Y' \rightarrow Y$ and $\alpha:X'\rightarrow X$ such that $Y'$ and $X'$
are nonsingular and
$f':Y'\rightarrow X'$ is a toroidal morphism, relative to
$\beta^{-1}(D_Y)$ and $\alpha^{-1}(D_X)$.
\end{them}

We need to generalize the notion of a monomial model defined in Definition 3 of
Section 2,
 to incorporate information about the divisors
$D_X$ and $D_Y$.

\begin{defin} A nonsingular model $f:Y\rightarrow X$ is said to be monomial
with respect to   SNC divisors $D_Y$
on $Y$ and $D_X$ on $X$ if $f^{-1}(D_X)=D_Y$ set theoretically, $C_f\subset
D_Y$, and if
for every $\beta \in Y$, with $\alpha := f(\beta) \in X$, there exist
regular systems of parameters  $(u,v)$ of
${\cal O}_{X,\alpha}$ and $(x,y)$ of ${\cal O}_{Y,\beta}$ such that
\begin{itemize}
\item[(i)]
If $\alpha \in {\rm Supp}(D_X)$, $D_X$ is locally at $\alpha$
defined by $u=0$ or $uv=0$.
\item[(ii)]
Either
$$
\left\{\begin{array}{lll}
u & = & \gamma x^a y^b \\
v & = & \delta x^c y^d
\end{array}
 \right. ,
$$
where $\gamma \delta$ is a unit in ${\cal O}_{Y,\beta}$, $p$ does not
divide $ad-bc$ and $D_Y$ is locally at $\beta$ defined by $xy=0$, or
$$
\left\{\begin{array}{lll}
u & = & \gamma x^a \\
v & = & \delta x^c \\
\end{array}  \right. ,
$$
where both of $\gamma \delta$ and
$a \gamma \displaystyle{\partial \delta \over \partial y} -
c \delta {\partial \gamma \over \partial y}$ are units in ${\cal O}_{Y,\beta}$, and $D_Y$ is locally at $\beta$ defined by $x=0$.
\end{itemize}
\end{defin}

\begin{them}
Suppose that $f:Y\rightarrow X$ is a  morphism of $k$-surfaces. Suppose that
$D_X$ and $D_Y$
are divisors on $X$ and $Y$ repectively, such that $f^{-1}(D_X)=D_Y$, and $D_Y$
contains all singular points of the
mapping $f$, and all singular points of $Y$. Then there exist projective
birational morphisms
 $\tau:Y_1\rightarrow Y$ and $\sigma:X_1\rightarrow X$
such that $Y_1$ and $X_1$ are nonsingular, the divisors
$D_{Y_1}=\tau^{-1}(D_Y)$ and $D_{X_1}=\sigma^{-1}(D_X)$
are SNC divisors, and $f_1:Y_1\rightarrow X_1$ is monomial with respect to
$D_{Y_1}$ and $D_{X_1}$.
\end{them}

The proof of Theorem 3  is a variation on the proof of the existence of a
monomial resolution, given in the preceeding sections.
Note that any model $f:Y\rightarrow X$ has divisors $D_X$ and $D_Y$ as in the
assumptions of the theorem.
After performing projective birational morphisms on $X$ and $Y$, we can assume
that
$X$ and $Y$ are nonsingular, and $D_X$ and $D_Y$ are SNC divisors. We must make
a change in the definition of
admissibility of Section 3.

Let $\alpha$ be a point in $X$. A r.s.p. $(u,v)$ of ${\cal O}_{X,\alpha}$
is said to be {\it admissible} if $\alpha \not \in {\rm Supp}(D_X)$,
or if $\alpha \in {\rm Supp}(D_X)$ and $D_X$ is locally at $\alpha$
defined by $u=0$ or $uv=0$ (see {\it (i)} in definition 6).\\

The results of Chapters 3 and 4 can now  be easily modified to produce a proof
of
Theorem 3.

By Theorem 3, we may
suppose that $X$, $Y$ are nonsingular $k$-surfaces, and that $f:Y\rightarrow X$
is  monomial
with respect to
 divisors $D_X$ on $X$ and $D_Y$ on $Y$. Thus  $D_X$ and $D_Y$
have simple normal crossings, $D_Y= f^{-1}(D_X)$ and
$C_f\subset D_Y$. We further have that for all $p\in D_X$
and $q\in f^{-1}(p)$ there exist regular parameters $(u,v)$ in ${\cal O}_{X,p}$
and $(x,y)$ in $\hat{\cal O}_{Y,q}$ such that one of the following holds.
\vskip .2truein
{\bf Case 1}
$D_{X,p} = V(uv)$, $D_{Y,q}=V(xy)$, $u=x^ay^b$, $v=x^cy^d$, $ad-bc\ne 0$.
\vskip .2truein
{\bf Case 2}
$D_{X,p} = V(uv)$, $D_{Y,q}=V(x)$, $u=x^a$, $v=x^c(y+\alpha)$, $0\ne \alpha\in k$.
\vskip .2truein
{\bf Case 3}
$D_{X,p} = V(u)$, $D_{Y,q}=V(xy)$, $u=x^ay^b$, $v=x^cy^d$, $ad-bc\ne 0$.
\vskip .2truein
{\bf Case 4}
$D_{X,p} = V(u)$, $D_{Y,q}=V(x)$, $u=x^a$, $v=x^c(y+\alpha)$, $\alpha\in k$.
\vskip .2truein
We will call cases 2 and 4 1-points, cases 1 and 3 2-points. Regular parameters
as above will be called
permissible.

The morphism is toroidal (relative to  $D_X$ and  $D_Y$) if all points satisfy
cases 1, 2 or 4*, where 4* is
\vskip .2truein
{\bf Case 4*}
$D_{X,p}=V(u)$, $D_{Y,q}=V(x)$, $u=x^a$, $v=y$.
\vskip .2truein
We will call a point $q\in Y$ good (or bad) if $f$ is toroidal (not toroidal)
at $q$.
By direct calculation, we see that

\begin{lem} The locus of bad points of $Y$ is closed of
 pure codimension 1 in $Y$.
\end{lem}

The set of image points in $X$ of bad points is finite.

Let
$$
G_f = \{q\in Y\vert q\in f^{-1}(p)\mbox{ is a 1-point such that }p\in X\mbox{
is the image of a bad point}\}.
$$
If $q\in G_f$, and $(x,y)$, $(u,v)$ are permissible parameters at $q$ and $p$,
we have an expression
$$
\begin{array}{l}
u=x^a\\
v=x^c(\gamma+y)
\end{array}
$$
for some $\gamma\in k$, and $D_{X,p}=V(u)$, $D_{Y,q}=V(x)$. We can define an
invariant for $q\in G_f$ by
$$
I(q,X)=\mbox{max}\{c-a\vert (x,y), (u,v)\mbox{ are permissible parameters at
$q$ and $p$}\}
$$
We then further define a global invariant
$$
r(Y,X) = \mbox{max}\{I(q,X)\vert q \in G_f\}.
$$

Suppose that $r(Y,X)>0$, and that $p\in X$ is such that there exists
$q\in f^{-1}(p)$ with $I(q,X)=r(Y,X)$.

Let $\pi:X_1\rightarrow X$ be the blowup of $p$. Let $f_1:Y\rightarrow X_1$
be the induced rational map, with $D_{X_1}=\pi^{-1}(D_X)$.

The following two lemmas are obtained by direct calculation of the effect of a
 quadratic transform at $q\in Y$ or at $p = f(q)\in X$.

\begin{lem} Suppose that $q\in f^{-1}(p)$ is such that $f_1$
is a morphism near $q$.
Suppose that $q$ is a 1-point. If $I(q,X)\le 0$, then $q$ is a good point for
$f_1$. If $I(q,X)>0$, then $I(q,X_1)<I(q,X)$.

The points where $f_1$ is not a morphism have one of the following forms.
$$
u=x^a, v=x^cy,\mbox{ with }c<a, (D_Y=V(x), D_X=V(u)), \mbox{ or}
\leqno (3)
$$
$$
u=x^ay^b, v=x^cy^d,\mbox{ with }a<c,b>d\mbox{ or }a>c, b<d,
 (D_Y=V(xy), D_X=V(u)).
\leqno (4)
$$
\end{lem}

\begin{lem} Suppose that $q\in f^{-1}(p)$ is such that $f_1$
is not a morphism near $q$.
Let $\tau:Y_1\rightarrow Y$ be the blowup of $q$, $f_2=f_1\circ\tau$.

Suppose that $q$ is a 1-point, so that $q$ has the form of (3).
Suppose that $q'\in \tau^{-1}(q)$. Then $f_2$ is a morphism at $q'$ and
$f_2$ is toroidal at $q'$ relative to $D_{Y_1}=\tau^{-1}(D_Y)$
and $D_{X_1}$, unless  $\hat{\cal O}_{Y_1,q'}$ has regular
parameters $x_1, y_1$ such that
$$
u=x_1^a, v=x_1^{c+1}y_1,
$$
with $D_{Y_1}=V(x_1), D_X=V(u)$, and
$$
I(q,X)<I(q',X)<0.
$$

Suppose that $q$ is a 2-point, so that $q$ has the form of (4).
We can assume that $a>c$ and  $b<d$.
Suppose that $q'\in \tau^{-1}(q)$.

If $q'$ is a 1-point,  then
$f_2$ is a morphism at $q'$ and either
$f_2$ is toroidal at $q'$ relative to $D_{Y_1}=\tau^{-1}(D_Y)$,
and $D_{X_1}$,
or we have $a+b<c+d$, so that  $\hat{\cal O}_{Y_1,q'}$ has regular
parameters $(x_1,  y_1)$, ${\cal O}_{X_1,f_2(q')}$
 has regular
parameters $(u_1,  v_1)$, such that
$$
u_1=u= x_1^{a+b},
 v_1=\frac{v}{u}= x_1^{c+d-(a+b)}(\gamma+ y_1),
$$
with $\gamma\ne 0$, $D_{Y_1}=V( x_1), D_{X_1}=V(u_1)$.
In this case we have
$$
I(q',X_1)=(c+d)-2(a+b)<(d-b)-1<r(Y,X).
$$
\end{lem}

By the above two lemmas, we can conclude

\begin{prop} There exist  sequences of quadratic transforms
$\alpha:Y' \rightarrow Y$ and $\beta:X'\rightarrow X$ such that
$f':Y'\rightarrow X'$ is a monomial mapping, and
$r(Y',X')\le 0$, relative to $\alpha^{-1}(D_Y)$ and $\beta^{-1}(D_X)$.
\end{prop}

Thus, we may assume that $r(Y,X)\le 0$. Suppose that $p\in X$ is the
image of a bad point. Then $v\mid u$ at all 2-points above $p$.
If $q$ is a 1-point above $p$, then
$$
u=x^a, v=x^c(\gamma+y)
$$
with $\gamma\in k$ and $c\le a$.
Let $\pi:X_1\rightarrow X$ be the blowup of $p$, $f_1:Y\rightarrow X_1$
be the induced rational map. Then $f_1$ is a morphism and is toroidal
at all 2-points of $f^{-1}(p)$, and at all 1-points with $\gamma\ne 0$,
and at all 1-points with $\gamma=0$ and  $c=a$.
The only points $q$ of $f^{-1}(p)$ where $f_1$ is not a morphism
(and is not toroidal) are 1-points of the form
$$
u=x^a, v=x^cy
$$
with  $I(q,X)=c-a<0$.
Let $\tau:Y_1\rightarrow Y$ be the blowup of such a $q$.
Let $f_2=f_1\circ \tau$. Then $f_2$ is a morphism and is toroidal at
all points of $\tau^{-1}(q)$ except possibly at a point $q'$
which has regular parameters $(x_1,y_1)$ satisfying
$x=x_1,y=x_1y_1$,
$$
u=x_1^a, v=x_1^{c+1}y_1
$$
with
$$
I(q,X)<I(q',X)<0
$$
By ascending induction on the negative number $I(q,X)$, we eventually
construct a toroidalization.
We thus have attained the conclusions of Theorem 2.

\vskip 5mm
\parindent=0mm
{\bf References}\\
\begin{itemize}\parindent=15mm
\item[\hbox to\parindent{\enskip\mbox{[Ab1]}}]
\textsc{Abhyankar S.}, On the valuations centered in a local domain,
{\it Amer. J. of Math.} {\bf 78} (1956), 321-348.
\item[\hbox to\parindent{\enskip\mbox{[Ab2]}}]
\textsc{Abhyankar S.}, Simultaneous resolution for algebraic
surfaces, {\it Amer. J. of Math.} {\bf 78} (1956), 761-790.
\item[\hbox to\parindent{\enskip\mbox{[Ab3]}}]
\textsc{Abhyankar S.}, Ramification theoretic methods in algebraic geometry,
{\it Annals of Math. Studies} {\bf 43},
Princeton University Press (1959).
\item[\hbox to\parindent{\enskip\mbox{[Ab4]}}]
\textsc{Abhyankar S.}, Local Uniformization on algebraic surfaces over ground
fields of characteristic $p\ne 0$,
{\it Annals  of Math.} {\bf 63} (1956), 491-526.
\item[\hbox to\parindent{\enskip\mbox{[Ab5]}}]
\textsc{Abhyankar S.}, Tame coverings and fundamental groups of algebraic varieties, 
{\it Amer. J. Math.} {\bf 81} (1959).
\item[\hbox to\parindent{\enskip\mbox{[AKa]}}]
\textsc{Abramovich D. and Karu K.},
Weak semistable reduction in characteristic 0, {\it preprint}.
\item[\hbox to\parindent{\enskip\mbox{[AKi]}}]
\textsc{Akbulut S. and King H.}, Topology of algebraic sets,
{\it MSRI publications} {\bf 25}, Springer Verlag (1992).
\item[\hbox to\parindent{\enskip\mbox{[C]}}]
\textsc{Cutkosky S.D.}, Local factorization and monomialization of morphisms,
{\it Ast{\'e}risque},  to appear.
\item[\hbox to\parindent{\enskip\mbox{[GM]}}]
\textsc{Grothendieck A. and Murre J.P.}, The tame fundamental group of a
formal neighbourhood of a divisor with normal crossings on a scheme,
{\it Lect. Notes in Math.} {\bf 208}, Springer-Verlag (1971).
\item[\hbox to\parindent{\enskip\mbox{[Ha]}}]
\textsc{Hartshorne R.}, Algebraic geometry, {\it Graduate Texts in Math.}
{\bf 52}, Springer-Verlag (1977).
\item[\hbox to\parindent{\enskip\mbox{[H]}}]
\textsc{Hironaka H.}, Desingularization of excellent surfaces, in
Cossart V., Giraud J. and Orbanz U., Resolution of singularities,
{\it Lect. Notes in Math.} {\bf 1101}, Springer-Verlag (1980).
\item[\hbox to\parindent{\enskip\mbox{[KKMS]}}]
\textsc{Kempf G., Knudsen F., Mumford D. and Saint-Donat B.},
Toroidal embeddings I, {\it Lect. Notes in Math.} {\bf 339},
Springer-Verlag (1973).
\item[\hbox to\parindent{\enskip\mbox{[Li]}}]
\textsc{Lipman J.}, Rational singularities, with applications to algebraic
surfaces and unique factorization,
{\it Publ. Math. IHES} {\bf 36} (1969), 195-279.
\item[\hbox to\parindent{\enskip\mbox{[Li2]}}]
\textsc{Lipman J.}, Desingularization of 2-dimensional schemes, {\it Annals
of Math.} {\bf 107} (1978), 115-207.
\item[\hbox to\parindent{\enskip\mbox{[Mi]}}]
\textsc{Milne J.S.}, Etale cohomology, {\it Princeton Series in Math.}
{\bf 33}, Princeton University Press (1980).
\item[\hbox to\parindent{\enskip\mbox{[SGA]}}]
\textsc{Grothendieck A.}, Rev\^etements \'etales et groupe fondemental, 
{\it Lecture Notes in Math} {\bf 224}, Springer Verlag, Heidelberg (1971).
\item[\hbox to\parindent{\enskip\mbox{[Sp]}}]
\textsc{Spivakovsky M.}, Valuations in function fields of surfaces,
{\it Amer. J. of Math.} {\bf 112} (1990), 107-156.
\item[\hbox to\parindent{\enskip\mbox{[Z1]}}]
\textsc{Zariski, O.}, Local uniformization of algebraic varieties, {\it Annals
of Math.},
 {\bf 41}, (1940), 852-896.
\item[\hbox to\parindent{\enskip\mbox{[ZS1]}}]
\textsc{Zariski O. and Samuel P.}, Commutive algebra I,
{\it Graduate Texts in Math.} {\bf 28}, Springer-Verlag (1958).
\item[\hbox to\parindent{\enskip\mbox{[ZS2]}}]
\textsc{Zariski O. and Samuel P.}, Commutive algebra II,
{\it The Univ.  Series in Higher Math.}, Van Norstrand (1960).
\par\end{itemize}

\end{document}